\documentclass{amsart}
\usepackage{amssymb,amsmath}

\newtheorem{theorem}{Theorem}
\newtheorem{lemma}{Lemma}
\newtheorem{example}{Example}
\newtheorem{problem}{Problem} 
\newcommand{\bt}{\begin{theorem}}
\newcommand{\et}{\end{theorem}}
\newcommand{\bl}{\begin{lemma}}
\newcommand{\el}{\end{lemma}}
\newcommand{\bex}{\begin{example}}
\newcommand{\eex}{\end{example}}
\newcommand{\bp}{\begin{problem}}
\newcommand{\ep}{\end{problem}}
\newcommand{\beal}{\begin{align*}}
\newcommand{\enal}{\end{align*}}
\newcommand{\beq}{\begin{equation}}
\newcommand{\eeq}{\end{equation}}
\newcommand{\benum}{\begin{enumerate}}
\newcommand{\eenum}{\end{enumerate}}
\newcommand{\ba}{\begin{array}}
\newcommand{\ea}{\end{array}}
\newcommand{\Z}{\ensuremath{\mathbf Z}}

\newcommand{\Q}{\ensuremath{\mathbf Q}}
\newcommand{\FF}{\ensuremath{\mathcal F}}

\newcommand{\mcs}{\ensuremath{\mathcal S}}
\newcommand{\mct}{\ensuremath{\mathcal T}}
\newcommand{\mcu}{\ensuremath{\mathcal U}}
\newcommand{\mcv}{\ensuremath{\mathcal V}}

\newcommand{\supp}{\text{supp}}

\newcommand{\pol}{$\mathcal{F} = \{f_n(q)\}_{n=1}^{\infty}$}
\newcommand{\polu}{$\mathcal{U} = \{u_n(q)\}_{n=1}^{\infty}$}
\newcommand{\polv}{$\mathcal{V} = \{v_n(q)\}_{n=1}^{\infty}$}

\begin{document}

\vspace{3cm}

\title{Linear quantum addition rules}
\author{Melvyn B. Nathanson}
\address{Department of Mathematics\\Lehman College (CUNY)\\Bronx, New York 10468}
\email{melvyn.nathanson@lehman.cuny.edu}
\dedicatory{To Ron Graham on his 70th birthday}
\thanks{This work was supported in part by grants from the NSA Mathematical Sciences Program
and the PSC-CUNY Research Award Program.}
\keywords{Quantum integers, quantum polynomial, quantum addition,
polynomial functional equation, $q$-series.}
\subjclass[2000]{Primary 11B37, 11P81, 65Q05, 81R50.  Secondary 11B13.}

\begin{abstract}
The quantum integer $[n]_q$ is the polynomial $1 + q + q^2 + \cdots + q^{n-1}.$
Two sequences of polynomials \polu\ and \polv\ define a {\em linear addition rule} $\oplus$ on 
a sequence \pol\ by $f_m(q)\oplus f_n(q) = u_n(q)f_m(q) + v_m(q)f_n(q).$
This is called a {\em quantum addition rule} if $[m]_q \oplus [n]_q = [m+n]_q$ for all positive integers $m$ and $n$.
In this paper all linear quantum addition rules are determined, and all solutions of the corresponding functional equations $f_m(q)\oplus f_n(q) = f_{m+n}(q)$ are computed. 
\end{abstract}

\maketitle

\section{Multiplication and addition of quantum integers}
We consider polynomials $f(q)$ with coefficients in a commutative ring with 1.
A sequence \pol\ of polynomials is nonzero if $f_n(q) \neq 0$ for some integer $n$.
For every positive integer $n$, the {\em quantum integer} $[n]_q$ is the polynomial
\[
[n]_q = 1 + q + q^2 + \cdots + q^{n-1}.
\]
These polynomials appear in many contexts.  In quantum calculus (Cheung-Kac~\cite{kac-cheu02}), for example, the $q$ derivative of $f(x) = x^n$ is
\[
f'(x) = \frac{f(qx) - f(x)}{qx-x} = [n]_q x^{n-1}.
\]
The quantum integers are ubiquitous in the study of quantum groups (Kassel~\cite{kass95}).

Let \pol\ be a sequence of polynomials.  Nathanson~\cite{nath03b}
observed that the multiplication rule
\[
f_m(q)\ast f_n(q) = f_m(q)f_n(q^m)
\]
induces a natural multiplication on the sequence of quantum integers,
since
\[
[m]_q \ast [n]_q = [mn]_q
\]
for all positive integers $m$ and $n$.  He asked what sequences \pol\ of polynomials, rational functions, and formal power series satisfy the multiplicative functional equation
\beq   \label{Qadd:femult}
f_m(q) \ast f_n(q) = f_{mn}(q)
\eeq
for all positive integers $m$ and $n$.
Borisov, Nathanson, and Wang~\cite{bori-nath-wang04} proved that the only solutions 
of~(\ref{Qadd:femult}) in the field $\Q(q)$ of rational functions with rational coefficients
are essentially quotients of products of quantum integers.  More precisely, let \pol\ be a nonzero 
solution of~(\ref{Qadd:femult}) in $\Q(q)$, and let $\supp(\FF)$ be the set of all integers $n$ with $f_n(q) \neq 0.$ 
They proved that there is a finite set $R$ of positive integers and a set $\{t_r\}_{r\in R}$ of integers such that, for all $n \in \supp(\mathcal{F})$,
\[
f_n(q) = \lambda(n) q^{t_0(n-1)}\prod_{r\in R} [n]_{q^r}^{t_r},
\]
where $\lambda(n)$ is a completely multiplicative arithmetic function and $t_0$ is a rational number 
such that $t_0(n-1) \in \Z$  for all $n \in \supp(\mathcal{F})$.
Nathanson~\cite{nath04c} also proved that if \pol\ is any solution of the functional
equation~(\ref{Qadd:femult}) in polynomials or formal power series with coefficients in a field, and if $f_n(0) = 1$ for all $n \in \supp(\FF),$ then there exists a formal power series $F(q)$ such that
\[
\lim_{n\rightarrow\infty \atop n \in \supp(\mathcal{F})} f_n(q) = F(q).
\]

Nathanson~\cite{nath03d} also defined the addition rule
\beq  \label{Qadd:addRule}
f_m(q) \oplus f_n(q) = f_m(q) + q^m f_n(q)
\eeq
on a sequence \pol\ of polynomials, 
and considered the additive functional equation
\beq  \label{Qadd:feAdd}
f_m(q) \oplus q^m f_n(q) = f_{m+n}(q).
\eeq
He noted that 
\beq   \label{Qadd:oplus}
[m]_q \oplus [n]_q  = [m+n]_q 
\eeq
for all positive integers $m$ and $n$,
and proved that every solution of the additive functional equation~(\ref{Qadd:feAdd}) is of the form
\[
f_n(q) = h(q)[n]_q,
\]
where $h(q) = f_1(q).$  This implies that if a nonzero sequence of polynomials \pol\ satisfies both the multiplicative functional equation~(\ref{Qadd:femult}) and additive function equation~(\ref{Qadd:feAdd}), then
\[
f_n(q) = [n]_q
\]
for all positive integers $n$.

In this paper we consider other binary operations $f_m(q)\oplus f_n(q)$ on sequences of polynomials that induce the natural addition of quantum integers or, equivalently, that satisfy~(\ref{Qadd:oplus}).
The goal of this paper is to prove that the addition rule~(\ref{Qadd:addRule}) is essentially the only linear quantum addition rule, and to find all solutions of the associated additive functional equation.

\section{Linear addition rules}
A general {\em linear quantum addition rule} is defined by two doubly infinite sequences of polynomials $\mcu = \{u_{m,n}(q)\}_{m,n=1}^{\infty}$ 
and $\mcv = \{v_{m,n}(q)\}_{m,n=1}^{\infty}$
such that
\beq   \label{Qadd:feGen}
[m+n]_q = u_{m,n}(q)[m]_q + v_{m,n}(q)[n]_q
\eeq
for all positive integers $m$ and $n$.
If the sequences $\mcu$ and $\mcv$ satisfy~(\ref{Qadd:feGen}), then $\mcu$ determines $\mcv$, and conversely.  
It is not known for what sequences $\mcu$ there exists a complementary sequence $\mcv$ satisfying~(\ref{Qadd:feGen}).

A {\em linear zero identity} is determined by two sequences of polynomials 
$\mcs = \{s_{m,n}(q)\}_{m,n=1}^{\infty}$ and $\mct = \{t_{m,n}(q)\}_{m,n=1}^{\infty}$ such that
\[
s_{m,n}(q)[m]_q + t_{m,n}(q)[n]_q = 0
\]
for all positive integers $m$ and $n$.

We can construct new addition rules from old rules  by adding zero identities and by taking affine combinations of addition rules.  
For example, the simplest quantum addition rule is
\beq  \label{Qadd:fe}
[m+n]_q = [m]_q + q^m[n]_q.
\eeq
Then
\[
[m]_q + q^m[n]_q = [m+n]_q = [n+m]_q = [n]_q + q^n[m]_q
\]
for all positive integers $m$ and $n$, 
and we obtain the zero identity
\beq   \label{Qadd:zero}
(1-q^n)[m]_q + (q^m-1)[n]_q = 0.
\eeq
Adding~(\ref{Qadd:fe}) and~(\ref{Qadd:zero}), we obtain
\beq  \label{Qadd:newrule}
[m+n]_q = (2-q^n)[m]_q + (2q^m-1)[n]_q.
\eeq
An affine combination of~(\ref{Qadd:fe}) and~(\ref{Qadd:newrule}) gives
\beq       \label{Qadd:example1}
[m+n]_q = (4-3q^n)[m]_q + (4q^m - 3)[n]_q.
\eeq

We can formally describe this process as follows.

\bt
For $i = 1,\ldots,k,$ let
$\mcu^{(i)} = \{u^{(i)}_{m,n}(q)\}_{m,n=1}^{\infty}$ and $\mcv^{(i)} = \{v^{(i)}_{m,n}(q)\}_{m,n=1}^{\infty}$ be sequences of polynomials that determine a quantum addition rule.
If $\alpha_1,\ldots, \alpha_k$ are elements of the coefficient ring such that 
\[
\alpha_1 + \cdots + \alpha_k = 1,
\]
and if the sequences $\mcu = \{u_{m,n}(q)\}_{m,n=1}^{\infty}$ and $\mcv = \{v_{m,n}(q)\}_{m,n=1}^{\infty}$ are defined by
\[
u_{m,n}(q) = \sum_{i=1}^k \alpha_i u^{(i)}_{m,n}(q)
\]
and
\[
v_{m,n}(q) = \sum_{i=1}^k \alpha_i v^{(i)}_{m,n}(q)
\]
for all positive integers $m$ and $n$, then the sequences $\mcu$ and $\mcv$ determine a quantum addition rule.

Similarly, if $\mcu = \{u_{m,n}(q)\}_{m,n=1}^{\infty}$ and $\mcv = \{v_{m,n}(q)\}_{m,n=1}^{\infty}$
are sequences of polynomials that determine a quantum addition rule, and if 
$\mcs = \{s_{m,n}(q)\}_{m,n=1}^{\infty}$ and $\mct = \{t_{m,n}(q)\}_{m,n=1}^{\infty}$ 
are sequences of polynomials that determine a zero identity, then the sequences $\mcu+S = \{u_{m,n}(q)+s_{m,n}(q)\}_{m,n=1}^{\infty}$ and $\mcv+T = \{v_{m,n}(q)+t_{m,n}(q)\}_{m,n=1}^{\infty}$
determine a quantum addition rule.
\et

\section{The fundamental quantum addition rule} 
In this paper we consider sequences $\mcu$ and $\mcv$ that depend only on $m$ or $n$.
We shall classify all linear zero identities and all linear quantum addition rules.

\bt      \label{Qadd:theorem:zerorule}
Let $\mcs = \{s_{n}(q)\}_{n=1}^{\infty}$ and $\mct = \{t_{m}(q)\}_{m=1}^{\infty}$ 
be sequences of polynomials.
Then
\beq    \label{Qadd:zerorule}
s_{n}(q)[m]_q + t_{m}(q)[n]_q = 0
\eeq
for all positive integers $m$ and $n$ if and only if there exists a polynomial $z(q)$ such that
\beq  \label{Qadd:zrs}
s_n(q) = z(q)[n]_q   \qquad\text{for all $n \geq 1$}
\eeq
and
\beq  \label{Qadd:zrt}
t_m(q) = -z(q)[m]_q   \qquad\text{for all $m \geq 1$.}
\eeq
If
\beq    \label{Qadd:zr3}
s_{m}(q)[m]_q + t_{m}(q)[n]_q = 0
\eeq
for all positive integers $m$ and $n$, 
or if
\beq   \label{Qadd:zr2}
s_{m}(q)[m]_q + t_{n}(q)[n]_q = 0
\eeq
for all positive integers $m$ and $n$, 
then $s_n(q) = t_n(q) = 0$ for all $n$.
\et

\begin{proof}
If there exists a polynomial $z(q)$ such that the sequences $\mathcal{S}$ and $\mathcal{T}$ satisfy identities~(\ref{Qadd:zrs}) and~(\ref{Qadd:zrt}), 
then
\[
s_n(q)[m]_q + t_m(q)[n]_q = z(q)[n]_q[m]_q-z(q)[m]_q[n]_q=0
\]
for all $m$ and $n$.

Conversely, suppose that the sequences \mcs\ and \mct\ define a linear zero identity of the form~(\ref{Qadd:zerorule}).
Letting $m = n = 1$ in~(\ref{Qadd:zerorule}), we have
\[
s_1(q) + t_1(q) = s_1(q)[1]_q + t_1(q)[1]_q = 0.
\]
Let
\[
z(q) = s_1(q) = -t_1(q).
\]
For all positive integers $n$ we have
\[
s_n(q)[1]_q + t_1(q)[n]_q =  s_n(q) - z(q)[n]_q = 0,
\]
and so
\[
s_n(q) = z(q)[n]_q.
\]
Similarly, 
\[
s_1(q)[m]_q + t_m(q)[1]_q = z(q)[m]_q + t_m(q)= 0,
\]
and so
\[
t_m(q) = -z(q)[m]_q
\]
for all positive integers $m$. 

If the sequences \mcs\ and \mct\ define a linear zero identity of the form~(\ref{Qadd:zr3}), then
\[
t_m(q)[n]_q = -s_m(q)[m]_q = t_m(q)[n+1]_q = t_m(q)([n]_q+q^n),
\]
and so 
\[
t_m(q)q^n = 0.
\]
It follows that $t_m(q)=0$ for all $m$, and so $s_m(q)=0$ for all $m$.

Suppose that the sequences \mcs\ and \mct\ define a linear zero identity of the form~(\ref{Qadd:zr2}).   Then 
\[
s_m(q)[m]_q = - t_n(q)[n]_q
\]
for all $m$ and $n$.  This implies that if $s_m(q) \neq 0$ for some $m$, then
 $t_n(q) \neq 0$ for all $n$ and $s_m(q) \neq 0$ for all $m$.  If \mcs\ and \mct\ are not the zero sequences, then, denoting the degree of a polynomial $f$ by $\deg(f),$ we obtain
\[
\deg(s_m)+m-1 = \deg(t_n)+n-1 \geq n-1,
\]
and so 
\[
\deg(s_m) \geq n-m
\]
for all positive integers $n$, which is absurd.  Therefore, \mcs\ and \mct\ are the zero sequences.
This completes the proof.
\end{proof}

\bt      \label{Qadd:theorem:mainiden}
Let $\mcu = \{u_{n}(q)\}_{n=1}^{\infty}$ and $\mcv = \{v_{m}(q)\}_{m=1}^{\infty}$ 
be sequences of polynomials.  Then
\beq      \label{Qadd:nm}
[m+n]_q = u_{n}(q)[m]_q + v_{m}(q)[n]_q
\eeq
for all positive integers $m$ and $n$ if and only if 
there exists a polynomial $z(q)$ such that
\beq   \label{Qadd:u}
u_n(q) = 1 + z(q)[n]_q  
\eeq
and
\beq   \label{Qadd:v}
v_m(q) = q^m - z(q)[m]_q 
\eeq
for all positive integers $m$ and $n$.
Moreover, $z(q) = u_1(q)-1 = q - v_1(q).$ 
\et

\begin{proof}
Let $z(q)$ be any polynomial, and define the sequences 
$\mcu = \{u_{n}(q)\}_{n=1}^{\infty}$ and $\mcv = \{v_{m}(q)\}_{m=1}^{\infty}$ 
by~(\ref{Qadd:u}) and~(\ref{Qadd:v}).  Then
\[
\begin{split}
u_n(q)[m]_q + v_m(q)[m]_q
& = \left( 1 + z(q)[n]_q\right) [m]_q + \left( q^m  - z(q)[m]_q \right) [n]_q \\
& = \left( [m]_q +  q^m [n]_q \right) + \left(  z(q)[n]_q [m]_q - z(q)[m]_q [n]_q \right) \\
& = [m]_q +  q^m [n]_q \\
& = [m+n]_q.
\end{split}
\]

Conversely, let $\mcu = \{u_{n}(q)\}_{n=1}^{\infty}$ and $\mcv = \{v_{m}(q)\}_{m=1}^{\infty}$ 
be a solution of~(\ref{Qadd:nm}).
We define 
\[
z(q) = u_1(q) - 1.
\]
Since
\[
1+q = [2]_q = [1+1]_q = u_1(q) + v_1(q) = 1 + z(q) + v_1(q),
\]
it follows that 
\[
v_1(q) = q - z(q).
\] 
For all positive integers $m$ we have
\[
[m+1]_q = u_1(q)[m]_q + v_m(q),
\]
and so 
\[
\begin{split}
v_m(q) & = [m+1]_q - u_1(q)[m]_q \\
& = q^m + [m]_q - u_1(q)[m]_q \\
& = q^m - z(q)[m]_q. 
\end{split}
\]
Similarly, for all positive integers $n$ we have
\[
[n+1]_q = [1+n]_q = u_n(q) + v_1(q)[n]_q,
\]
and so
\[
\begin{split}
u_n(q) & = [n+1]_q - v_1(q)[n]_q \\
& = 1 + q[n]_q - (q - z(q))[n]_q \\
& = 1 + z(q)[n]_q.
\end{split}
\]
This completes the proof.
\end{proof}

For example, we can rewrite the quantum addition rule~(\ref{Qadd:example1}) in the form
\[
\begin{split}
[m+n]_q & = (4-3q^n)[m]_q + (4q^m - 3)[n]_q  \\
& = (1 + z(q)[n]_q)[m]_q + (q^m - z(q)[m]_q)[n]_q,
\end{split}
\]
where
\[
z(q) = 3-3q.
\]

\bt           \label{Qadd:theorem:secondiden}
Let $\mcu = \{u_{m}(q)\}_{m=1}^{\infty}$ and $\mcv = \{v_{n}(q)\}_{n=1}^{\infty}$ 
be sequences of polynomials.
Then
\beq    \label{Qadd:qar-mm}
[m+n]_q = u_{m}(q)[m]_q + v_{m}(q)[n]_q
\eeq
for all positive integers $m$ and $n$ if and only if $u_m(q) = 1$ and $v_m(q) = q^m$ for all $m$.
There do not exist sequences of polynomials $\mcu = \{u_{m}(q)\}_{m=1}^{\infty}$ and $\mcv = \{v_{n}(q)\}_{n=1}^{\infty}$ such that
\beq        \label{Qadd:qar-mn}
[m+n]_q = u_{m}(q)[m]_q + v_{n}(q)[n]_q
\eeq
for all positive integers $m$ and $n$.
\et

\begin{proof}
Suppose that for every positive integer $m$ we have
\[
[m+1]_q = u_m(q)[m]_q + v_m(q)[1]_q = u_m(q)[m]_q + v_m(q),
\]
and
\[
[m+2]_q = u_m(q)[m]_q + v_m(q)[2]_q = u_m(q)[m]_q + (1+q)v_m(q).
\]
Subtracting, we obtain
\[
q^{m+1} = [m+2]_q - [m+1]_q = q v_m(q),
\]
and so
\[
v_m(q) = q^m.
\]
Then
\[
u_m(q)[m]_q = [m+1]_q - v_m(q) = [m+1]_q - q^m = [m]_q,
\]
and so
\[
u_m(q) = 1
\]
for all $m$.
This proves the first assertion of the Theorem.

If~(\ref{Qadd:qar-mn}) holds for $n = 1$ and all $m$, then
\[
[m+1]_q = u_m(q)[m]_q + v_1(q)[1]_q = u_m(q)[m]_q + v_1(q),
\]
and so
\[
u_m(q)[m]_q = [m+1]_q - v_1(q).
\]
We also have
\[
[m+2]_q = u_m(q)[m]_q + v_2(q)[2]_q = [m+1]_q - v_1(q) + (1+q)v_2(q),
\]
and so
\[
q^{m+1} = [m+2]_q - [m+1]_q = (1+q)v_2(q) - v_1(q)
\]
for all positive integers $m$, which is absurd.
\end{proof}

Theorems~\ref{Qadd:theorem:mainiden} and~\ref{Qadd:theorem:secondiden} show that all linear quantum addition rules are of the form $[m+n]_q = u_{n}(q)[m]_q + v_{m}(q)[n]_q.$  The following result shows that the sequence of quantum integers is essentially the only solution of the corresponding functional equation.

\bt            \label{Qadd:theorem:AddFE}
Let $\mcu = \{u_{n}(q)\}_{n=1}^{\infty}$ and $\mcv = \{v_{m}(q)\}_{m=1}^{\infty}$ 
be sequences of polynomials such that
\[
[m+n]_q = u_{n}(q)[m]_q + v_{m}(q)[n]_q
\]
for all positive integers $m$ and $n$.
Then \pol\ is a solution of the functional equation
\[
f_{m+n}(q)= u_{n}(q)f_m(q) + v_{m}(q)f_n(q)
\]
if and only if there is a polynomial $h(q)$ such that 
\[
f_n(q) = h(q) [n]_q 
\]
for all $n \geq 1.$
\et

\begin{proof}
By Theorem~\ref{Qadd:theorem:mainiden},
there exists a polynomial $z(q)$ such that
\[
u_n(q) = 1 + z(q)[n]_q  
\]
and
\[
v_m(q) = q^m - z(q)[m]_q 
\]
for all positive integers $m$ and $n$.
The proof is by induction on $n$.
Let $h(q) = f_1(q).$
Suppose that $f_n(q) = h(q)[n]_q$ for some integer $n \geq 1.$
Then
\[
\begin{split}
f_{n+1}(q) 
& = u_1(q) f_n(q) + v_n(q) f_1(q) \\
& = (1 + z(q) ) h(q)[n]_q + ( q^n - z(q)[n]_q ) h(q) \\
& = h(q)([n]_q+q^n)\\
& =h(q) [n+1]_q . 
\end{split}
\]
This completes the proof.
\end{proof}

{\em Remark.}
The only property of polynomials used in this paper is the degree of a polynomial, which occurs in the proof that there is no nontrivial zero identity of the form~(\ref{Qadd:zr2}).  It follows that Theorems~\ref{Qadd:theorem:mainiden}  and~\ref{Qadd:theorem:AddFE}  hold in any algebra that contains the quantum integers, for example, the polynomials, the rational functions, the formal power series, or the formal Laurent series with coefficients in a ring or field.  

\section{Nonlinear addition rules}
A. V. Kontorovich observed that the quantum integers satisfy the following two nonlinear addition rules:
\[
[m+n]_q =  [m]_q + [n]_q - (1-q)[m]_q[n]_q
\] 
and
\[
[m+n]_q = q^n [m]_q + q^m [n]_q + (1-q)[m]_q[n]_q.
\]
These give rise to the functional equations
\[
f_m(q)\oplus f_n(q) =  f_m(q) + f_n(q) - (1-q)f_m(q)f_n(q)
\]
and
\[
f_m(q)\oplus f_n(q) = q^n f_m(q) + q^mf_n(q) + (1-q)f_m(q)f_n(q),
\]
whose solutions are, respectively,
\[
\begin{split}
f_n(q) &  = \frac{1}{q-1}\sum_{k=1}^n {n\choose k} ((q-1)f_1(q))^k\\
& = \frac{ 1 - (1+(q-1)f_1(q))^n }{1 - q}.
\end{split}
\]
and
\[
\begin{split}
f_n(q) &  = \frac{1}{q-1}\sum_{k=1}^n {n\choose k} q^{n-k} ((1-q)f_1(q))^k\\
& = \frac{ (q+(1-q)f_1(q))^n - q^n }{1-q}.
\end{split}
\]
Kontorovich and Nathanson~\cite{kont-nath05} have recently described all quadratic addition rules for the quantum integers.
It would be interesting to classify higher order nonlinear quantum addition rules.

\providecommand{\bysame}{\leavevmode\hbox to3em{\hrulefill}\thinspace}
\providecommand{\MR}{\relax\ifhmode\unskip\space\fi MR }
\providecommand{\MRhref}[2]{%
  \href{http://www.ams.org/mathscinet-getitem?mr=#1}{#2}
}
\providecommand{\href}[2]{#2}

\end{document}